\documentclass[12 pt]{amsart}
\title[La gravedad, ?`Una propiedad esencial...]{La Gravedad, ?`Una propiedad esencial de la materia? La controversia Newton-Leibnitz.}
\author{Jonathan Taborda}
\email{taborda50@gmail.com}
\date{Mayo de 2010}
\usepackage{amsmath,amssymb,amsthm,pb-diagram,mathpazo}
\usepackage[activeacute,spanish]{babel}
\usepackage[T1]{fontenc}
\begin{document}\maketitle
%\begin{abstract}
%\end{abstract}
\fontfamily{ppl}\selectfont
El \textit{Philosophiae naturalis principia mathematica}, más popularmente conocido con el nombre de \textit{Principia Mathematica}, es sin duda alguna la obra cumbre y la síntesis de la Revolución Científica de los siglos XVI-XVII iniciada con la publicación del \textit{De Revolutionibus} por el Polaco N. Copernicus; publicada por vez primera en London en 1687, por el Lucasian Profesor Sir Isaac Newton (1643-1727). Este voluminoso texto, consta de una estructura tripartita, i.e, tres libros que inicialmente fueron redactados por Newton como Lecciones Universitarias en la University of Cambridge. El primero libro cubre el tópico del movimiento de los cuerpos, el segundo sobre el movimiento de los cuerpos en un medio resistente y el tercero y más polémico intitulado El sistema del Mundo, a ocasionado todo tipo de reacciones y/o polémicas desde su primera edición. La causa de dicho litigio radica, como es bien sabido en la no explicación mecanicista por parte del Sabio Inglés a la fuerza de atracción Gravitacional en el libro III del \textit{Principia}.\\
Una de las múltiples controversias que se generaron a lo largo de la $1^a$ Ed., la constituye la lectura y/o interpretación que del \textit{Principia} realizara el polifacético y gran Matemático G.W.Von Leibnitz,\footnote{Cf. \textit{La controversia Newton-Leibnitz.} Karen Gloy. Carlos Emel Rendón-Traductor. Estudios de Filosofía. Febrero de 1993.} quién tiene conocimiento de la obra Newtoniana durante su estancia en París (1672-1676).A decir verdad, fue Leibnitz quien, a pesar de su propio uso del concepto de atracción o de su juego con él, en su \textit{Tentamen de motuum coelestium causis} de 1689, o quizá debido a ello, acentuó la analogía entre la atracción y una cualidad oculta, y la reforzó afirmando que la atracción era un milagro. Lo hizo, empero, bastante tarde.\\
Al principio, en 1690\footnote{En una carta a Huygens, Octubre de 1690. Cf. Alexandre Koyré. \textit{Estudios Newtonianos}. Traducción, introducción y notas por Felipe Ochoa. Univ. de Antioquia. Inst. de Filosofía. 1998. Apéndice B. pp. 310. ss.} sólo expresó su sorpresa de que Newton pareciera considerar la gravedad como una virtud incorpórea e inexplicable y no creyera en la posibilidad de explicarla mecánicamente-?`No lo había hecho Huygens muy bien?.\\
Más en 1690 Leibnitz no habló respecto de cualidades ocultas o acerca de milagros. De hecho, fue Huygens a quien reprochó por introducir milagros en la filosofía natural: ?`No es la indivisibilidad de los átomos, en la que Huygens creía, un milagro perpetuo? Sólo en 1703, en su \textit{Nouveaux essays} (que, no obstante, no publicó debido a la muerte de Locke: no quería atacar aun adversario fallecido) menciona ''cualidades ocultas'' y ''milagros''.  De este modo, al referirse a la carta de Locke a Stillingflect, Leibnitz escribe:
\begin{quote}
\textit{<<A continuación cita el ejemplo de la atracción de la materia (pp. 99, pero sobre todo pp. 408.), y habla de la gravitación de la materia hacia la materia, reconociendo que nunca se podrá concebir cómo se da la gravitación. Lo cual es tanto como volver a las cualidades ocultas, o inclusive más, inexplicables>>.}\footnote{Cf. Koyré. Op. Cit.}
\end{quote}
Sólo en 1711, en una carta a N. Hartsoeker,\footnote{Febrero 6 de 1711, publicadas primero en las \textit{Mémoires de Trévoux} (1711), y luego en una traducción inglesa, en las \textit{Memoirs of Literature} de 1712. Cf. Ibid.} hace esto.\\
En 1715 (noviembre o diciembre), en una carta al Abbé Conti, Leibnitz renovó su ataque:
\begin{quote}
\textit{<<Su filosofía me parece bastante extraña y no puedo creer que se pueda justificar. Si todo cuerpo es pesado, se sigue (digan lo que digan sus partidarios, aunque lo nieguen apasionadamente) que la gravedad será una cualidad oculta escolástica, o es más, el efecto de un milagro...No es suficiente con decir: Dios ha hecho una ley de la naturaleza, en consecuencia la cosa es natural. Es necesario que la ley la pueda explicar la naturaleza de las cosas creadas. Si por ejemplo, Dios fuera a dar a un cuerpo libre la ley de girar alrededor de cierto centro, tendría que  juntar bien fuera éste cuerpo con otros los cuales por su impulso lo harían permanecer siempre en una órbita circular, o ponerlo en los talones de un ángel; o más, tendría que coincidir extraordinariamente con su movimiento. Pues naturalmente saldría por la tangente... Estoy fuertemente a favor de la filosofía experimental, pero el Sr. Newton se aleja mucho de ella cuando pretende que toda la materia sea pesada (o que cada parte de la materia atraiga a otra) lo cual ciertamente no está probado por la experimentación...\\
Y puesto que no conocemos todavía en detalle cómo se produce la gravedad, o la fuerza elástica, o fuerza magnética, ello no nos da derecho alguno para hacer de ellas cualidades ocultas, escolásticas, o milagros; incluso tampoco nos da derecho a limitar la sabiduría de Dios.>>}\footnote{Cf. Ibid. pp. 314}\\
\end{quote}
Por ese mismo tiempo, el problema de la atracción, i.e, la cuestión de sí es una  cualidad oculta y un milagro, o una respetable fuerza de la naturaleza, formó uno de los temas principales de la famosa polémica entre Leibnitz y Clarke, junto con el problema de la realidad o imposibilidad del espacio vacío, el movimiento absoluto y otros problemas de metafísica y filosofía natural.\\
Es bien sabido que Newton no creyó que la gravedad fuera una <<propiedad innata, esencial e inherente de la materia>>. De hecho, en 1675, en su \textit{Hypothesis Explaining the properties of light} y en 1679, en una carta a Boyle (28 de febrero de 1678/9), intentó explicar la gravedad por medios mecánicos- i.e, por los movimientos de la materia sutil o un medio etéreo- pero no prosiguió estos vanos intentos, al menos por algún tiempo; en 1692, en una carta a Bentley, le pidió no atribuirle esa noción epicúrea; la atracción como acción a distancia a través del vacío sin mediación, decía a Bentley, era un gran absurdo en que nadie podría creer y además, afirmó de forma muy clara, que esta medición tenía que ser ejecutada por algo que no es material, o sea, por Dios.\footnote{Cf. Ibid.}\\
La acción a distancia se ha considerado muy amenudo como un medio de explicación inaceptable en la física.\\
Las acciones a distancia siempre fueron una característica de la tradición magica natural a través de la Edad Media. En distintos pasajes de sus escritos, Isaac Newton da por hecho la posibilidad de la acción a distancia. No parece tener dificultad alguna con el punto de vista según el cual un cuerpo puede afectar a otro aun cuando estén separados en el espacio. Considérense, V,g., las palabras iniciales de las muy influyentes Cuestiones con que Newton concluye su Óptica:<<?`No actúan los cuerpos sobre la luz a \textit{distancia}, y por su acción flexionan sus rayos>>? (Cuestión
1).\footnote{Cf.\textit{Isaac Newton y el problema de la acción a distancia} John Henry. Traduccion:Felipe Ochoa. Revista Estudios de Filos. Febrero de 2007. Univ. de Antioquia. pp. 189-226.} Así pues, grandes eruditos del pensamiento Newtoniano como lo son Rupert Hall, Alexandre Koyré e I. Bernhard Cohen, afirman que Newton no creía que la gravedad pudiera ser una propiedad inherente de la materia, ni que la acción a distancia fuera posible <<sin mediacion>>. Para Epicuro la atracción gravitacional era una propiedad esencial de la materia; para Newton esto era inconcebible porque una propiedad tal requería un Mediador inmaterial.\footnote{Cf. John Henry. Op. cit. pp. 119}\nocite{*}\\
<<Para Newton y Bentley Dios era el mediador inmaterial cuya omnipotencia le permitia imponerle a la materia un agente secundario de atracción gravitacional que actúa constantemente de acuerdo a ciertas leyes>>.
La influencia de Newton no se limita al campo matemático y físico, sino que se extiende, como puede obviamente deducirse además al campo filosófico y también al teológico.\\
Como hemos subrayado anteriormente, los puntos de vista matemáticos, físicos, filosóficos y teológicos de Newton, además de estar condensados en el \textit{Principia}, también están embebidos en su famosa \textit{Optics, or A treatise of the Reflexions, Reflections, Inflections and Colors of the light} (1704, ampliados por 31 Cuestiones).\footnote{Cf. Gloy. Op. cit. pp. 10 y ss.} En el \textit{Principia Mathematica} las concepciones newtonianas relativas a la absolutez del espacio-tiempo, independiente de las cosas que ocupan el espacio y de la materia, la <<caja del mundo>>, generaron la amarga polémica Leibnitz-Newton, ó siendo más precisos Leibnitz-Clarke, ya que éste discípulo de Newton y editor de la edición Latina de la Óptica, es quién se encarga de la controversia filosófica con el cortesano Hannoveriano.\\
Lo característico del espacio Newtoniano es su status ontológico i.e., la suposición de su ser real. Por eso Newton afirma:
\begin{quote}
<<\textit{El espacio no es afección de un cuerpo o de otro cuerpo, de cualquier ser finito, ni pasa de sujeto a sujeto, sino que es siempre invariablemente la inmensidad de un sólo y siempre el mismo inmensum. Los espacios finitos no son en modo alguno afecciones de sustancias finitas, sino que son solamente las partes del espacio infinito en las que existen las sustancias finitas.\\
El espacio no está limitado por los cuerpos, sino que existe igualmente tanto con cuerpos como sin ellos. El espacio no está encerrado entre cuerpos, sino que éstos, existiendo en un espacio ilimitado, están solamente delimitados por sus propias dimensiones.\\
El espacio infinito es uno, absoluto y esencialmente indivisible, y suponerlo dividido es una contradicción en los términos, porque en su separación debe haber espacio, lo cual es suponerlo dividido, y sin embargo, no dividido al mismo tiempo.\\
El espacio es uniforme o igual y no difiere una parte de otra.\\
Además, el espacio y el tiempo son cantidades.}>>\footnote{Cf. Ibid.}

\end{quote}
A la concepción del sabio Inglés se opone, la del mayor Spinozista de todos los tiempos; para Leibnitz, espacio y tiempo son sistemas materiales de relación con un correspondiente carácter relativo. Por tanto, éste define el espacio como << orden de la existencia en la simultaneidad>> y el tiempo como <<orden de la sucesión.>>\footnote{Cf. Ibid. pp. 13} Leibnitz concibe el espacio como un sistema de posiciones, lo cual corresponde a su ideal de ciencia del Análisis.\\
La consideración de un espacio absoluto en Newton, radicará en que sólo un espacio absoluto puede explicar el fenómeno de la simetría de derecha e izquierda, arriba y abajo, adelante atrás, como en general, el fenómeno de la existencia separada de objetos congruentes.\\
A tales argumentos, Leibnitz opone un principio lógico-metafísico, conocido como el de razón suficiente, en el que concidera que Dios <<no había tenido motivo y/o razón suficiente para crear dos veces el mismo mundo.>>\\
<<Si Newton ha conquistado méritos por la elaboración de un estracto fenomenológico de la espacialidad, Leibnitz los ha alcanzado por el tratamiento matemático conjuntual del espacio.>>\\

A continuación, nos enfocaremos un poco en el \textit{General Scholium} insertado por Newton a la segunda Ed. del \textit{Principia} (1713), porque allí es donde Newton hace manifiesta y explícita su adhesión a las doctrinas Sociano-Antitrinitarias y responde el porque de su no creencia en la gravedad como una <<propiedad innata, esencial e inherente de la materia>>.\\
En una carta bien conocida (del 18 de febrero, 1712/3), Cotes llamó la atención de Newton hacia el ataque que Leibnitz hizo en su contra, publicado en la \textit{Memoirs of Literature} de Mayo de 1712, y le aconsejó no dejarlo sin responder.\footnote{Cf. Koyré. Op. cit. Cap. IV. pp. 451. ss.} Asimismo, Cotes formuló una objeción a la teoría de Newton sobre la atracción o, al menos, sobre la manera como fue presentada. Seguramente Cotes era consciente de las numerosas y frecuentes afirmaciones de que Newton usaba el término <<atracción>> como término perfectamente neutro, <<conjuntamente>> con otros, y que se podía entender como si significara presión o cualquier otra cosa, salvo lo que parecía significar. No obstante, Cotes notó que la atracción newtoniana implicaba la atribución de <<fuerzas de atracción>> a los cuerpos, y que Newton, tácitamente, hacía esa <<hipótesis>> o <<suposición>>.\\
Estas objeciones por parte del esbirro Cotes al sabio Inglés, desenbocarán en el prefacio que éste compa\~{n}ero de armas tomar en la <<Hermandad Polaca>> o comunidad Sociana redactara en 1713, con el objetivo de responder a los críticos continentales.\\
El III libro del \textit{principia}: El sistema del mundo, inicia con un conjunto de Reglas para el estudio de la Filosofía Natural. Es de vital importancia, subrayar el hecho de que sus famosas \textit{regulae philosophandi} que aparecen en dicho texto, son una versión transformada, de un conjunto de 16 en total que aparecen en el \textit{Treatise on the Apocalypse} y que en el \textit{Principia} de Newton únicamente dos del conjunto originario en el \textit{Treatise} aparecen en la $1^a$ Ed. de 1687, donde ellas son parte de la <<hipótesis>> introductoria.\\
Newton sentó nueve <<hipótesis>> como base para su <<Sistema del Mundo>> en el libro III. En la segunda Ed. (1713) denominó <<phenomena>> a las cinco últimas hipótesis consistentes en afirmaciones acerca de los planetas y de la luna, y a las tres primeras las llamó <<\textit{Regulae Philosophandi}>>, que se puede traducir como <<Reglas de inferencias en las ciencias naturales>>.\footnote{Cf. \textit{De las leyes de Kepler a la ley de Grav. Usal. de Newton y viceversa}. Jaime Chica, con la asistencia de Jonathan Taborda. Cap. I. Metamorfosis de la Ciencia. pp. 20. ss. en preparación.}\\
Luego de las mencionadas 4 reglas en la $2^{da}$ Ed. del \textit{Principia}, aparece el celebérrimo \textit{Escolio General}.\\
El Escolio General del \textit{Principia Mathematica} ha sido
caracterizado como "posiblemente el más famoso de todos los escritos
Newtonianos", como lo afirma el Prof. I.B.Cohen en su
\textit{Introduction to Newton's Principia}. Sí esto es así, el
Escolio General es la porción mejor-conocida de uno de los trabajos
más importantes en la Historia de la Ciencia. No cabe duda que el
Escolio contiene tres de las líneas más citadas de Newton. Aquí se
proclama el "más bello Sistema para el Sol, Planetas y Cometas,
únicamente procede del consejo y dominio para un Ser inteligente y
poderoso". El Escolio General también ofrece el pronunciamiento para
el discurso de Dios; "de la apariencia de las cosas, ciertamente se
hacen pertenecer a la Filosofía Natural", a lo largo del
desconocimiento Newtoniano en la afirmación sobre la causa de la
gravedad "frame no hypotheses" \textit{hypotheses non fingo}.\\
Dos estudios recientes han avanzado considerablemente sobre nuestro
conocimiento para la dimensión posterior del Escolio General.
Primero, en un ensayo sobre el Dios del Dominio de Newton, James
Force punteó la naturaleza antitrinitaria \index{antitrinitarismo} contenida en éste.
Segundo, Larry Stewart presentó una narración detallada para quien
el Escolio fue leído y conocido como un documento heterodoxo por el
más perceptivo de los enemigos de los amigos de Newton. Stewart
también construye un fuerte caso para mostrar que Newton estuvo
ciertamente empleando su Escolio General para soportar públicamente
a su aliado Samuel Clarke, cuya \textit{Scripture-Doctrine of the
Trinity} no-ortodoxa había aparecido en 1712\footnotemark
\footnotetext{Cf. "God of gods, and Lord of Lords". The Theology of Isaac
Newton General Scholium to the Principia.\\
Stephen D. Snobelen. The History of Science Society. 2001.}.\\
Una de las presentaciones más sorprendentes para el Dios del Escolio
General es su pensamiento Hebraico y carácter bíblico. Newton no
duda que su Dios no fue otro que el "Dios de Israel". Todos los
ejemplos de nombres divinos, títulos y atributos dados en el Escolio
son citados directamente, ó alusiones no ambiguas para las
Escrituras. Algunas de las expresiones, tales como "Señor Dios",
''el Dios de Israel", "mi Dios'' y "nuestro Dios", ocurren en la
Escritura también como una lista detallada. Otros ejemplos son de
especial consideración. El título ''Señor Dios $\pi\alpha\nu\tau
o\kappa\rho\alpha\tau\omega\varrho$'' aparece en seis ocasiones en
el Nuevo Testamento, únicamente en el libro de la Revelación- un
libro que tuvo una fascinación particular para la mentalidad
profética de Newton\footnotemark\footnotetext{Los seis ejemplos para
éste título están en Rev. 4.8, 11.17, 15.3, 16.7, 19.6 y 21.22. Cf.
Ibid. pp.
177. n. 34.}.\\
El título "Señor de señores" (cuando lo aplicó para el Padre),
aparece en el Salmo 136.3 y Timoteo 6.15. Sólo que Newton es aún más
específico. La presentación de Dios en el Escolio es estrictamente
unitaria y monoteísta en el sentido Hebreo.\\
Arguyendo que "Dios'' es un término definido por relaciones, Newton
rechaza específicamente las definiciones absolutas de Dios, como
''eterno'', ''infinito'' y ''perfecto''.\\
Él ciertamente acepta esos términos como atributos de Dios y así lo
afirma en el Escolio.\\
Excepto que él también deja en claro que ellos son inadecuados como
sinónimos de "Dios", porque ellos no expresan relaciones. Además,
tales palabras son problemáticas como definiciones fijas para el
término ''Dios'' en una segunda dirección, entonces en el sentido
estricto, ellos no son grados de eternidad, infinitud y perfección,
y por tanto no se pueden aplicar para un ser que no sea otro que el
Dios verdadero.\footnote{Cf. Jaime Chica, Jonathan Taborda, Op. cit.}\\
Consideramos pues, que es demasiado pristino y tajante el \textit{hypotheses non fingo}, plasmado en el Escolio General por parte del Sabio Inglés, para responder a las críticas y malas interpretaciones que de la concepción Gravitacional se han hecho a lo largo de la historia del \textit{Principia}. También estamos deacuerdo con la idea del Prof. Henry consistente en afirmar que <<Betty Jo Dobbs también se convenció, de que Newton no pudo haber creído en las acciones a distancia, y por tanto sugirió que las Cuestiones de la Óptica en las que se refirió a las acciones a distancia debieron haber sido diseñadas deliberadamente para ocultar sus especulaciones privadas reales, en las cuales según Dobbs implicaban una creencia en el espíritu vegetal de la materia.>>\footnote{Cf.Henry. Op. cit. pp. 116.}La Profesora Dobbs, sostiene que la gravedad Newtoniana se derivaba de ideas neo-estoicas en las que el neuma estoico se platoniza en un principio inmaterial-el cual Newton, de paso, igualó a Dios. <<Lo que llamamos gravitación-escribe ella- no es nada más que la Mente Divina formando y reformando las partes del Universo>>.\footnote{Cf. Henry. Op. cit. n. 29.} Henry en su documento, se opone tajantemente a la interpretación canónica de los expertos Hall, Cohen, Koyré y Dobbs; ya que para ellos <<la gravedad es llevada a cabo por la potentia absoluta de Dios, su poder absoluto.>> Para Henry <<la gravedad está a cargo de la potentia ordinaria de Dios, su poder ordenado.>>\\ En su interpretación, Newton negó explícitamente que Dios fuera la causa de la gravedad; <<En él están contenidas y se mueven todas las cosas, pero él no actúa sobre ellas ni ellas sobre Él. Dios no experimenta nada por los movimientos de los cuerpos; los cuerpos no sienten resistencia por la omnipresencia de Dios.>>\footnote{Cf. Ibid. pp. 204.} Para Newton, pues <<la gravitación es una energía constante infundida a la materia por el Autor de todas las cosas.>>
En el Libro III (1713), Newton deduce, a partir de la reinterpretación y correción\footnote{Cf. Ibid. Cap. I. pp. 6. n. 17.} de las tres leyes del movimiento planetario postuladas por Kepler en \textit{Astronomia Nova} (1609) y \textit{Harmonices Mundi} (1619) una descripción teórica para una fuerza \underline{matemática}, más no física, conocida como Ley de Gravitación Universal.
\bibliographystyle{amsalpha}
\bibliography{arias}

\end{document}